\documentclass[graybox]{svmult}

\usepackage{helvet}         % selects Helvetica as sans-serif font
\usepackage{courier}        % selects Courier as typewriter font
\usepackage{type1cm}        % activate if the above 3 fonts are
\usepackage{makeidx}         % allows index generation
\usepackage{graphicx}        % standard LaTeX graphics tool
\usepackage{multicol}        % used for the two-column index
\usepackage[bottom]{footmisc}% places footnotes at page bottom

\usepackage{placeins}
\usepackage[cp1251]{inputenc}
\begin{document}

\title*{The Quantum Mirror to the Quartic del Pezzo Surface}
\author{H\"ulya Arg\"uz}
\institute{H\"ulya Arg\"uz \at Institute of Science and Technology, Austria, \email{huelya.arguez@ist.ac.at}
%\and Name of Second Author \at Name, Address of Institute %\email{name@email.address}
}
\maketitle

\abstract{
A log Calabi--Yau surface $(X,D)$ is given by a smooth projective surface $X$, together with an anti-canonical cycle of rational curves $D \subset X$. The homogeneous coordinate ring of
the mirror to such a surface -- or to the complement $X\setminus D$ -- is constructed using wall structures and is generated by theta functions \cite{GHK, GHS}. In \cite{A1}, we provide a recipe to concretely compute these theta functions from a combinatorially constructed wall structure in $\vec{R}^2$, called the heart of the canonical wall structure. In this paper, we first apply this recipe to obtain the mirror to the quartic del Pezzo surface, denoted by $dP_4$, together with an anti-canonical cycle of $4$ rational curves. We afterwards describe the deformation quantization of this coordinate ring, following \cite{BP}. This gives a non-commutative algebra, generated by quantum theta functions. There is a totally different approach to construct deformation quantizations using the realization of 
the mirror as the monodromy manifold of the Painlev\'e IV equation \cite{CMR,M}. We show that these two approaches agree.
}

\section{The mirror to $(dP_4,D)$}
\label{sec:2}
\begin{figure}[b]
\sidecaption
% For example, with the graphicx style use
\includegraphics[scale=.3]{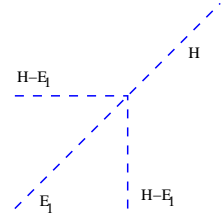}
%
% If no graphics program available, insert a blank space i.e. use
%\picplace{5cm}{2cm} % Give the correct figure height and width in cm
%
\caption{The kinks of the MVPL function $\varphi$ on the rays of the fan $\Sigma_{\overline{X}}$ associated to the toric model $\overline{X}$, drawn in dashed blue lines to distinguish them in what follows from walls.}
\label{fig:PL}     % Give a unique label
\end{figure}

We construct the mirror to $(dP_4,D)$, following the recipe in \cite{A1}. To obtain the theta functions generating the coordinate ring of the mirror, we define an initial wall structure associated to $(X,D)$, using the following data:

\begin{itemize}
    \item A choice of a \emph{toric model}, that is, a birational morphism $(X, D) \to (\overline{X},\overline{D})$ to a smooth
toric surface $\overline{X}$ with its toric boundary $\bar{D}$ such that $D \to \overline{D}$ is an isomorphism. For the quartic del Pezzo surface $X=dP_4$, obtained by blowing up $5$ general points in ${\vec{P}}^2$, we consider the toric model given by a toric blow-up of ${\vec{P}}^2$, illustrated as in Figure \ref{Fig: dP40}.

We choose $\overline{D}$ to be the toric boundary divisor in $\overline{X}$, and let $D$ be the strict transform of $\overline{D}$. We note that the equations of the mirror will be independent of the choice of the toric model \cite{A1,GHK}. 

This choice defines a natural subdivision of $M_{\vec{R}}$, where $M \cong \vec{Z}^2$ is a fixed lattice, and $M_{\vec{R}} = M \otimes\vec{R} \cong \vec{R}^2$, given by the toric fan $\Sigma_{\overline{X}} \subset M_{\vec{R}}$ of $\overline{X}$. We denote $M_{\vec{R}}$ together with this subdivision by $(M_{\vec{R}},\Sigma_{\overline{X}})$.

\begin{figure}[b]
\includegraphics[scale=.3]{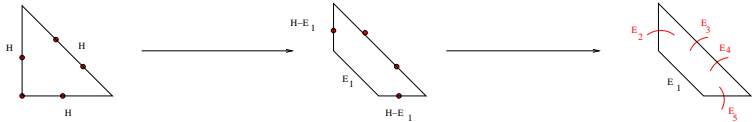}
\caption{On the left hand figure $H$ is the class of a general line on the projective plane $\vec{P}^2$, which is illustrated by its momentum map image. In the middle $\overline{X}$ is the the blow up of $\vec{P}^2$ at a toric point, $E_1$ is the class of the exceptional curve, and by abuse of notation $H$ denotes the class of the pull-back of a general line. On the right hand figure we illustrate $dP_4$ together with the exceptional divisors obtained by $4$ further non-toric blow-ups.}
\label{Fig: dP40}      % Give a unique label
\end{figure}

  \item A \emph{multi-valued
piecewise linear} (MVPL) function $\varphi$ on $(M_{\vec{R}},\Sigma_{\overline{X}})$ with values in the monoid of integral points of the cone of effective curves on $X$, denoted by $\mathrm{NE}(X)$. Up to a linear function, we uniquely define $\varphi$ by specifying its kinks along each ray of $\Sigma_{\overline{X}}$ to be the pullback of the class of the curve in $\overline{X}$ corresponding to this ray, as illustrated in Figure \ref{fig:PL}.
\end{itemize}

%\begin{minipage}{0.5\textwidth}

%\end{minipage} \hfill
%\begin{minipage}{0.65\textwidth}

%\begin{figure}[H]
%\center{\resizebox{.4\linewidth}{!}{\input{kinks.pspdftex}}}
%\caption{\label{fig:PL} The kinks of the MVPL function}
%\end{figure}
%\end{minipage}

\begin{definition}
\label{def wall structure}
A wall structure on $(M_{\vec{R}},\Sigma_{\overline{X}})$ is a collection of pairs $(\rho,f_{\rho})$, called walls, consisting of rays $\rho \subset M_{\vec{R}}$, together with functions $ f_{\rho} \in \vec{C}[NE(X)^{\mathrm{gp}}][M]$, referred to as wall-crossing functions. Each wall crossing function defines a wall-crossing transformation 
\begin{equation}
\label{Eqn: crossing walls}
\theta_{\gamma,\rho}: \quad
z^v\longmapsto f_\rho^{\langle n_\rho,v \rangle} z^v \,,
\end{equation}
prescribing how a monomial $z^v$ changes when crossing the wall $(\rho,f_{\rho})$. Here, $n_\rho$ is the normal vector to $\rho$ chosen with a sign convention as in \cite[\S$2.2.1$]{AG}. 
\end{definition}

We note that in general one might need to work with infinitely many walls, and then %to keep track of the powers of $t$, 
to work modulo an ideal of $NE(X)$ while defining the wall crossing functions. However, in the case of the quartic del Pezzo surface we will only need finitely many walls, and the wall crossing functions will be elements of a polynomial ring, as in Definition \ref{def wall structure}. 

To obtain the equation of the mirror to $(X,D)$, we first define an initial wall structure associated to $(X,D)$. To do this, we first define an initial set of walls in  $(M_{\vec{R}},\Sigma_{\overline{X}})$. 
For every non-toric blow-up in the toric model, we include a wall 
$(\rho,f_\rho)$, where $\rho$ is the ray in $\Sigma_{\overline{X}}$ 
corresponding to the divisor on which we do the non-toric blow-up, 
and
\[ f_{\rho} = 1+ t^{-E_i}z^{-v_{\rho}} = 1+ t^{-E_i}x^{-a}y^{-b} \,,\]
where $E_i \in NE(X)$ is the class of the exceptional curve, $z^{v_{\rho}}$ denotes the element in the monoid ring $\vec{C}[M]$, corresponding to the primitive direction $v_{\rho} \in M$ of $\rho$ pointing towards $0$ referred to as the direction vector, and we use the convention $x= z^{(1,0)}$ and $y= z^{(0,1)}$. The negative signs on the powers of $t$ and $z$ are chosen following the sign conventions of \cite{AG}.

\begin{figure}[b]
\includegraphics[scale=.25]{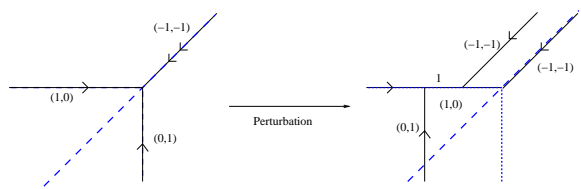}
\caption{Walls of the initial wall structure associated to $(dP_4,D)$, and their perturbations obtained by translating some walls, so that each intersection is formed by only two walls meeting at a point. The dashed blue lines indicate where the kinks of the MVPL function $\varphi$ are, the black rays with arrows are walls. On the left hand figure the function attached to wall with direction $(1,0)$ is given by $1+t^{-E_2}x^{-1}$, the function on the wall wit direction $(0,1)$ is $1+t^{-E_5}y^{-1}$, and the one on the wall with direction $(-1,-1)$ is $(1+t^{-E_3}xy)(1+t^{-E_4}xy)$, as there are two walls on top of each other. On the right hand figure the functions on the walls with $(1,0)$ and $(0,1)$ remain same. However, we have now two separate walls in direction $(-1,-1)$, with attached functions $(1+t^{-E_3}xy)$ and $(1+t^{-E_4}xy)$ respectively -- we choose to attach the latter to the most right ray. 
%(the choice of which one to attach to which will not change the resulting computation of the mirror).
}
\label{Fig: perturb}     % Give a unique label
\end{figure}

The walls of the initial wall structure do not generally intersect only pairwise, but there can be triple or more complicated intersections as illustrated on the left hand Figure \ref{Fig: perturb}, where we have $4$ walls intersecting - two of these walls lie on top of each other and have direction vector $(-1,-1)$, and the other two have direction vectors $(1,0)$ and $(0,1)$.  However, we can always move these walls, so that any of the intersection points of the initial walls will be formed by only two walls intersecting. The initial wall structure after a choice of such a perturbation is illustrated on the right hand side of Figure \ref{Fig: perturb}. In this situation, only when two walls intersect each time, we can easily describe a \emph{consistent wall structure} obtained from the initial wall structure by inserting new walls to it, so that the composition of the wall crossing transformations around each intersection point is identity. This diagram is referred as the heart of the canonical wall structure in \cite{A1}.

 Roughly put, each time a wall with support on a ray $\rho_i$ with direction $v_i$ and with attached function $ 1+t^{-E_i}z^{-v_i}$ intersects another wall with support on a ray $\rho_j$ with direction $v_j$ and with attached function $1+t^{-E_j}z^{-v_j}$, we first extend these rays to lines, so that on the new rays we form while doing this we attach the functions $1+t^{\kappa_i-E_i}z^{-v_i}$ and $1+t^{\kappa_j-E_j}z^{-v_j}$ respectively. Here $\kappa_i$ and $\kappa_j$ denote the sums of the kinks of the MVPL function $\varphi$, lying on rays that intersect the initial rays. We furthermore insert an additional wall with support on a ray $\rho_{i+j}$ with direction $v_i+v_j$ and with attached function \[ 1+t^{(\kappa_i+\kappa_j)-(E_i+E_j)}z^{-(v_i+v_j)} \,.\]
Note that this simple prescription describes the consistent wall structure, because each time two walls with direction vectors say $v_i$ and $v_j$ intersect, the determinant of $v_i$ and $v_i$ is $\pm 1$. 
%%Note that in the particular example considered in this paper we are in situation where we have only finitely many walls in the consistent wall structure, with all wall crossings simple, therefore it suffices to 

\begin{figure}[b]
\sidecaption
\includegraphics[scale=.2]{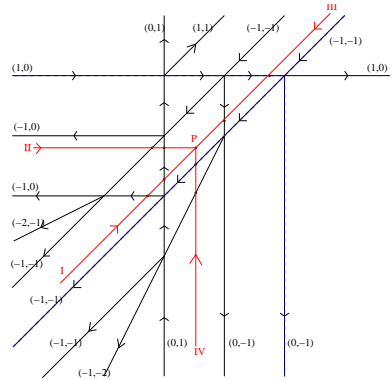}
\caption{The consistent wall structure obtained from the perturbed initial wall structure with $4$ incoming rays, illustrated in Figure \ref{Fig: perturb}. The kinks of the MVPL function $\varphi$ are on the $4$ blued dashed rays, whereas all the black rays are walls. Tracing the red rays labelled by I, II, III, and IV, and the walls crossed by them while they reach the point $P$, we determine the theta $4$ functions $\vartheta_1,\vartheta_2,\vartheta_3,\vartheta_4$.}  
\label{Fig: DGPScons}
% Give a unique label
\end{figure}

The coordinate ring of the mirror to $(X,D)$ is generated by \emph{theta functions}, determined by keeping track of how a set of initial monomials change under wall crossing in this consistent wall structure \cite[\S3]{A1}. 
%We illustrate this construction in the case of $(dP_4,D)$. First, note that 
We have as many theta functions as the number of rays of the toric model associated to $(dP_4,D)$, whose fan is illustrated in Figure \ref{fig:PL}. The direction vectors of these rays are given by $(-1,-1),(-1,0),(1,1),(0,-1)$, which correspond to the initial set of monomials $x^{-1}y^{-1},x^{-1},xy,y^{-1}$. The corresponding $4$ theta functions, respectively denoted by $\vartheta_1,\ldots,\vartheta_4$, are determined by tracing how these monomials change as the corresponding rays cross walls in the consistent wall structure associated to $(dP_4,D)$ illustrated in Figure \ref{Fig: DGPScons}. To determine the wall crossings we first fix a general point $P$ as in Figure \ref{Fig: DGPScons} 
%(the construction of the mirror will be independent of this choice)
, and look at the rays coming from the directions $(-1,-1),(-1,0),(1,1),(0,-1)$ and stopping at this point -- note that, we made a choice of the point $P$ so that each ray will cross exactly two walls and the situation will be symmetric, and it will make it easier to calculate the theta functions. 

Each of the $4$ theta functions is then obtained by the following wall crossings:  The red ray labelled with I crosses the two walls with attached functions $1+t^{H-E_4-E_5}x$ and $1+t^{E_1-E_5}y^{-1}$, the functions on the two walls crossed by the red ray labelled by II are $1+t^{E_1-E_5}y^{-1}$ and $1+t^{H-E_1-E_3}xy$, the functions on the walls crossed by the red ray labelled by III are $1+t^{-E_2}x^{-1}$ (in this case there is also the kink of the MVPL function $H-E_1$ we take into account in the computation of the theta functions) and $1+t^{H-E_1-E_2-E_3}y$, and finally the red ray labelled by IV crosses the walls with attached functions $1+t^{2H-E_1-E_2-E_3-E_4}xy^2$ and $1+t^{H-E_1-E_4}xy$ (and in this case there is also a kink of the MVPL function given by $E_1$). Hence, we obtain;
\begin{eqnarray}
\nonumber
& &     x^{-1}y^{-1}  \mapsto    x^{-1}y^{-1} +t^{H-E_4-E_5} y^{-1} \mapsto x^{-1}y^{-1} +t^{E_1-E_5} x^{-1} y^{-2} + t^{H-E_4-E_5} y^{-1}   =:   \vartheta_1  \\
\nonumber
    & &   x^{-1}   \mapsto    x^{-1} +t^{E_1-E_5} x^{-1} y^{-1} \mapsto x^{-1} +t^{H-E_1-E_3} y + t^{E_1-E_5} x^{-1}y^{-1}   =:  \vartheta_2 \\
     \nonumber
     & &   xy   \mapsto    t^{H-E_1} xy + t^{H-E_1-E_2} y  \mapsto t^{H-E_1} xy + t^{2H-2E_1-E_2-E_3} xy^2 +t^{H-E_1-E_2} y  =:   \vartheta_3 \\
     \nonumber
& &         y^{-1}   \mapsto    y^{-1} +t^{2H-E_1-E_2-E_3-E_4} xy \mapsto t^{E_1}y^{-1} +t^{H-E_4} x + t^{2H-E_1-E_2-E_3-E_4} xy   =:   \vartheta_4 
\nonumber
\end{eqnarray}

The above theta functions generate the mirror to  $(dP_4,D)$, which is obtained as a family of log Calabi--Yau surfaces $(dP_4,D)$ given
as 
$\mathrm{Spec}$ of the quotient of 
$\vec{C}[\mathrm{NE}(X)][\vartheta_1,\vartheta_2,\vartheta_3,\vartheta_4]$
by the following two quadratic equations:
\begin{eqnarray}
\label{Eq: classic quadric eqns}
 \vartheta_1 \vartheta_3 & = & C_1 + t^{H-E_1-E_2} \vartheta_2 + t^{H-E_1-E_5} \vartheta_4 \\
 \nonumber
 \vartheta_2 \vartheta_4 & = & C_2 + t^{E_1} \vartheta_1 + t^{H-E_3-E_4} \vartheta_3 
 \nonumber
\end{eqnarray} 
where
\begin{eqnarray}
\label{Eq C1C2}
C_1 & = & t^{H-E_1} + t^{2H-E_1-E_2-E_3-E_5} +  t^{2H-E_1-E_2-E_4-E_5}  \\
\nonumber
C_2 & = & t^{H-E_4} + t^{H-E_3} +  t^{2H-E_2-E_3-E_4-E_5}.  
\end{eqnarray} 

Note that the resulting equations for the mirror given above, agrees with the one obtained in \cite{B} using computer algebra. 
%However, here we provide the concrete description of each of the theta functions, which will be useful while determining the deformation quantization in the following section.
%By doing this computation by hand here we are able to provide the concrete description of the wall structure as well as each of the theta functions, which will be useful while determining the deformation quantization in the following section.
 
\section{The quantum mirror to $(dP_4,D)$}
A general recipe to construct a deformation quantization of mirrors of log Calabi--Yau surfaces is given in \cite{BP}. For $(dP_4,D)$, as the consistent wall structure consists of only finitely many walls, the general recipe reduces to the following simple prescription:
the quantum theta functions, denoted by $\hat{\vartheta}_1,\ldots,\hat{\vartheta}_4$, are obtained from the theta functions above by replacing the monomials $z^v \in \vec{C}[M]$, by quantum variables, denoted by  $\hat{z}^v$, which are elements of the quantum torus, that is such that 
$\hat{z}^v \hat{z}^{v'}=q^{\frac{1}{2}\det(v,v')}\hat{z}^{v+v'}$, where $q$ is the quantum deformation parameter. Hence, from the equations for the theta functions above, we obtain:
\begin{eqnarray}
\hat{\vartheta}_1 & = & \hat{z}^{(-1,-1)} + t^{E_1-E_5} \hat{z}^{(-1,-2) }  + t^{H- E_4-E_5} \hat{z}^{(0,-1)} \\
\nonumber
\hat{\vartheta}_2 & = & \hat{z}^{(-1,0)} + t^{H-E_1-E_3} \hat{z}^{(0,1) }  + t^{ E_1-E_5} \hat{z}^{(-1,-1)} \\
\nonumber
\hat{\vartheta}_3 & = & t^{H-E_1}\hat{z}^{(1,1)} + t^{2H-2E_1-E_2-E_3} \hat{z}^{(1,2) }  + t^{H- E_1-E_2} \hat{z}^{(0,1)} \\
\nonumber
\hat{\vartheta}_4 & = & t^{E_1} \hat{z}^{(0,-1)} + t^{H-E_4} \hat{z}^{(1,0) }  + t^{2H-E_1-E_2-E_3-E_4} \hat{z}^{(1,1)} 
\end{eqnarray}

These quantum theta functions satisfy the followsing $8$ relations -- first two obtained as deformations of the two quadric equations in (\ref{Eq: classic quadric eqns}), and the latter $6$ equations determining the non-commutativity of the products of each two among four quantum theta functions.

\begin{eqnarray}
\label{Eq:quantul theta}
\hat{\vartheta}_1 \hat{\vartheta}_3 & = & C_1 + q^{-\frac{1}{2}} t^{H-E_1-E_2} \hat{\vartheta}_2 +  q^{\frac{1}{2}} t^{H-E_1-E_5} \hat{\vartheta}_4 \\
\nonumber
\hat{\vartheta}_2 \hat{\vartheta}_4 & = & C_2 + q^{\frac{1}{2}} t^{E_1} \hat{\vartheta}_1 +  q^{-\frac{1}{2}} t^{H-E_3-E_4} \hat{\vartheta}_3 \\
\nonumber
q^{\frac{1}{2}} \hat{\vartheta}_1 \hat{\vartheta}_3 - q^{-\frac{1}{2}} \hat{\vartheta}_3 \hat{\vartheta}_1 & = & (q^{\frac{1}{2}}- q^{-\frac{1}{2}})C_1+(q-q^{-1})t^{H-E_1-E_5}\hat{\vartheta}_4 \\
\nonumber
q^{\frac{1}{2}} \hat{\vartheta}_2 \hat{\vartheta}_4 - q^{-\frac{1}{2}} \hat{\vartheta}_4 \hat{\vartheta}_2 & = & (q^{\frac{1}{2}}- q^{-\frac{1}{2}})C_2+(q-q^{-1})t^{E_1}\hat{\vartheta}_1 \\
\nonumber
q^{\frac{1}{2}} \hat{\vartheta}_1 \hat{\vartheta}_2 - q^{-\frac{1}{2}} \hat{\vartheta}_2 \hat{\vartheta}_1 & = & (q^{\frac{1}{2}}- q^{-\frac{1}{2}}) t ^{2H-E_1-E_3-E_4-E_5} \\
\nonumber
q^{\frac{1}{2}} \hat{\vartheta}_2 \hat{\vartheta}_3 - q^{-\frac{1}{2}} \hat{\vartheta}_3 \hat{\vartheta}_2 & = & (q^{\frac{1}{2}}- q^{-\frac{1}{2}}) t ^{H-E_5} \\
\nonumber
q^{\frac{1}{2}} \hat{\vartheta}_3 \hat{\vartheta}_4 - q^{-\frac{1}{2}} \hat{\vartheta}_4 \hat{\vartheta}_3 & = & (q^{\frac{1}{2}}- q^{-\frac{1}{2}}) t ^{H-E_2} \\
\nonumber
q^{\frac{1}{2}} \hat{\vartheta}_4 \hat{\vartheta}_1 - q^{-\frac{1}{2}} \hat{\vartheta}_1 \hat{\vartheta}_4 & = & (q^{\frac{1}{2}}- q^{-\frac{1}{2}}) t ^{2H-E_1-E_2-E_3-E_4} 
\end{eqnarray}

Eliminating $\hat{\vartheta}_4$ in (\ref{Eq:quantul theta}), we obtain the same equations as in \cite[Corollary 4.6]{M}, obtained there using a totally different approach based on the fact that the cubic equation obtained from (\ref{Eq: classic quadric eqns}) by eliminating $\vartheta_4$ is the defining equation of the wild character variety arising as the monodromy manifold of the Painlev\'e IV equation.

\FloatBarrier

\end{document}